\documentclass[12pt]{article}
\usepackage{amssymb,amsmath,amsthm,tikz,graphicx}
\usetikzlibrary{arrows}

\title{On Deformed Dodecahedron Tiling}
\author{Yohji Akama, 
Tohoku University \\
Min Yan\thanks{Research was supported by Hong Kong RGC General Research Fund 605610 and 606311.},
Hong Kong University of Science and Technology}

\theoremstyle{definition}

\theoremstyle{remark}

\begin{document}

\maketitle

\begin{abstract}
There is only one type of tilings of the sphere by $12$ congruent pentagons. These tilings are isohedral.
\end{abstract}

A deformed dodecahedron tiling is a tiling of the sphere by $12$ congruent pentagons. Gao, Shi, and Yan completely classified these tilings into five types in \cite{gsy}, with the fifth type allowing two free parameters. Since the first four types can only contain isolated examples, it was conjectured in \cite{gsy} that the four types actually can only allow the regular dodecahedron tiling. In this note, we confirm the conjecture, so that the fifth type in \cite{gsy} is exactly all the deformed dodecahedron tilings. We will also show that the symmetry of any deformed dodecahedron tiling is transitive on the tiles.

\bigskip

\begin{figure}[htp]
\centering
\begin{tikzpicture}

\draw
	(0,1) node[below=1] {\small $\alpha$} 
		-- (-1,0) node[right] {\small $\beta$}
		-- (-0.7,-1) node[above right=-1] {\small $\delta$}
		-- (0.7,-1) node[above left=-1] {\small $\epsilon$}
		-- (1,0) node[left] {\small $\gamma$} 
		-- (0,1)
	(0.4,0.4) -- (0.6,0.6)
	(0.45,0.35) -- (0.65,0.55)
	(-0.4,0.4) -- (-0.6,0.6)
	(-0.45,0.35) -- (-0.65,0.55)
	(0.7,-0.4) -- (1,-0.5)
	(-0.7,-0.4) -- (-1,-0.5);

\end{tikzpicture}
\caption{Geometrical constraint for pentagon.}
\label{pent}
\end{figure}
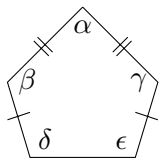

\noindent{\bf 1.} 
Cheuk, Cheung, and Yan give a geometrical constraint for spherical pentagons in \cite[Lemma 3]{ccy1}. Their lemma says that, if the spherical pentagon in Figure \ref{pent} has two pairs of equal edges, then $\beta>\gamma$ if and only if $\delta<\epsilon$.

\begin{figure}[htp]
\centering
\begin{tikzpicture}[scale=1]

\foreach \a in {0,1}
{
\begin{scope}[xshift=7*\a cm]

\foreach \x in {1,...,5}
{
\draw[rotate=72*\x]
	(90:0.8) -- (18:0.8) -- (18:1.5) -- (54:2) -- (90:1.5)
	(54:2) -- (54:3) -- (126:3);
	
\coordinate (X\x A\a) at (18+72*\x:0.5);
\coordinate (Y\x A\a) at (6+72*\x:0.95);
\coordinate (Z\x A\a) at (-38+72*\x:0.95);
\coordinate (U\x A\a) at (8+72*\x:1.35);
\coordinate (V\x A\a) at (-44+72*\x:1.35);
\coordinate (W\x A\a) at (-18+72*\x:1.65);
\coordinate (P\x A\a) at (18+72*\x:1.7);
\coordinate (Q\x A\a) at (-11+72*\x:2.1);
\coordinate (R\x A\a) at (47+72*\x:2.1);
\coordinate (S\x A\a) at (-13+72*\x:2.6);
\coordinate (T\x A\a) at (49+72*\x:2.6);
\coordinate (O\x A\a) at (-18+72*\x:3.2);
\coordinate (P\a) at (-60:3.5);
}

\end{scope}
}



\node at (X1A0) {\small $\alpha$};
\node at (X2A0) {\small $\beta$};
\node at (X3A0) {\small $\delta$};
\node at (X4A0) {\small $\epsilon$};
\node at (X5A0) {\small $\gamma$};

\node at (Y1A0) {\small $\gamma$};
\node at (Z1A0) {\small $\alpha$};
\node at (U1A0) {\small $\epsilon$};
\node at (V1A0) {\small $\beta$};
\node at (W1A0) {\small $\delta$};

\node at (Y2A0) {\small $\beta$};
\node at (Z2A0) {\small $\delta$};
\node at (U2A0) {\small $\alpha$};
\node at (V2A0) {\small $\epsilon$};
\node at (W2A0) {\small $\gamma$};

\node at (Y3A0) {\small $\epsilon$};
\node at (Z3A0) {\small $\gamma$};
\node at (U3A0) {\small $\delta$};
\node at (V3A0) {\small $\alpha$};
\node at (W3A0) {\small $\beta$};

\node at (Y4A0) {\small $\delta$};
\node at (Z4A0) {\small $\epsilon$};
\node at (U4A0) {\small $\beta$};
\node at (V4A0) {\small $\gamma$};
\node at (W4A0) {\small $\alpha$};

\node at (Y5A0) {\small $\delta$};
\node at (Z5A0) {\small $\epsilon$};
\node at (U5A0) {\small $\beta$};
\node at (V5A0) {\small $\gamma$};
\node at (W5A0) {\small $\alpha$};

\node at (P1A0) {\small $\delta$};
\node at (Q1A0) {\small $\epsilon$};
\node at (R1A0) {\small $\beta$};
\node at (S1A0) {\small $\gamma$};
\node at (T1A0) {\small $\alpha$};

\node at (P2A0) {\small $\alpha$};
\node at (Q2A0) {\small $\beta$};
\node at (R2A0) {\small $\gamma$};
\node at (S2A0) {\small $\delta$};
\node at (T2A0) {\small $\epsilon$};

\node at (P3A0) {\small $\alpha$};
\node at (Q3A0) {\small $\beta$};
\node at (R3A0) {\small $\gamma$};
\node at (S3A0) {\small $\delta$};
\node at (T3A0) {\small $\epsilon$};

\node at (P4A0) {\small $\beta$};
\node at (Q4A0) {\small $\delta$};
\node at (R4A0) {\small $\alpha$};
\node at (S4A0) {\small $\epsilon$};
\node at (T4A0) {\small $\gamma$};

\node at (P5A0) {\small $\gamma$};
\node at (Q5A0) {\small $\alpha$};
\node at (R5A0) {\small $\epsilon$};
\node at (S5A0) {\small $\beta$};
\node at (T5A0) {\small $\delta$};

\node at (O1A0) {\small $\alpha$};
\node at (O2A0) {\small $\gamma$};
\node at (O3A0) {\small $\epsilon$};
\node at (O4A0) {\small $\delta$};
\node at (O5A0) {\small $\beta$};


\draw[line width=1.5, xshift=7cm]
	(-54:0.8) -- (234:0.8)
	(18:0.8) -- (18:1.5)
	(162:0.8) -- (162:1.5)
	(-18:2) -- (-18:3)
	(198:2) -- (198:3)
	(54:3) -- (126:3);


\node at (X1A1) {\small $\alpha$};
\node at (X2A1) {\small $\beta$};
\node at (X3A1) {\small $\delta$};
\node at (X4A1) {\small $\epsilon$};
\node at (X5A1) {\small $\gamma$};

\node at (Y1A1) {\small $\beta$};
\node at (Z1A1) {\small $\delta$};
\node at (U1A1) {\small $\alpha$};
\node at (V1A1) {\small $\epsilon$};
\node at (W1A1) {\small $\gamma$};

\node at (Y2A1) {\small $\epsilon$};
\node at (Z2A1) {\small $\gamma$};
\node at (U2A1) {\small $\delta$};
\node at (V2A1) {\small $\alpha$};
\node at (W2A1) {\small $\beta$};

\node at (Y3A1) {\small $\gamma$};
\node at (Z3A1) {\small $\epsilon$};
\node at (U3A1) {\small $\alpha$};
\node at (V3A1) {\small $\delta$};
\node at (W3A1) {\small $\beta$};

\node at (Y4A1) {\small $\epsilon$};
\node at (Z4A1) {\small $\delta$};
\node at (U4A1) {\small $\gamma$};
\node at (V4A1) {\small $\beta$};
\node at (W4A1) {\small $\alpha$};

\node at (Y5A1) {\small $\delta$};
\node at (Z5A1) {\small $\beta$};
\node at (U5A1) {\small $\epsilon$};
\node at (V5A1) {\small $\alpha$};
\node at (W5A1) {\small $\gamma$};

\node at (P1A1) {\small $\alpha$};
\node at (Q1A1) {\small $\beta$};
\node at (R1A1) {\small $\gamma$};
\node at (S1A1) {\small $\delta$};
\node at (T1A1) {\small $\epsilon$};

\node at (P2A1) {\small $\gamma$};
\node at (Q2A1) {\small $\alpha$};
\node at (R2A1) {\small $\epsilon$};
\node at (S2A1) {\small $\beta$};
\node at (T2A1) {\small $\delta$};

\node at (P3A1) {\small $\gamma$};
\node at (Q3A1) {\small $\epsilon$};
\node at (R3A1) {\small $\alpha$};
\node at (S3A1) {\small $\delta$};
\node at (T3A1) {\small $\beta$};

\node at (P4A1) {\small $\beta$};
\node at (Q4A1) {\small $\alpha$};
\node at (R4A1) {\small $\delta$};
\node at (S4A1) {\small $\gamma$};
\node at (T4A1) {\small $\epsilon$};

\node at (P5A1) {\small $\beta$};
\node at (Q5A1) {\small $\delta$};
\node at (R5A1) {\small $\alpha$};
\node at (S5A1) {\small $\epsilon$};
\node at (T5A1) {\small $\gamma$};

\node at (O1A1) {\small $\delta$};
\node at (O2A1) {\small $\epsilon$};
\node at (O3A1) {\small $\gamma$};
\node at (O4A1) {\small $\alpha$};
\node at (O5A1) {\small $\beta$};

\end{tikzpicture}
\caption{Types 1 and type 4 tilings.}
\label{T14}
\end{figure}
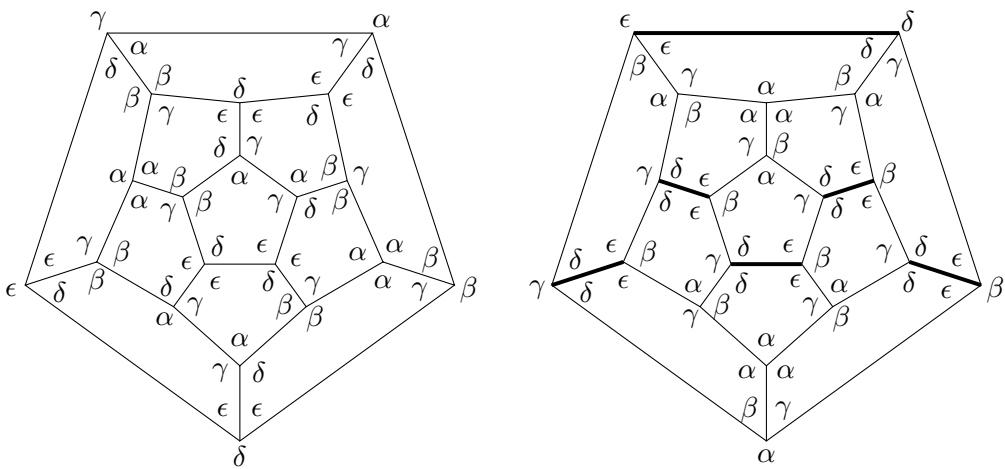

The type 1 tiling of \cite{gsy} is given on the left of Figure \ref{T14} (up to permutation of angle notations), with all edges having equal length. If $\epsilon>\alpha$, then by the fact that the sum of angles at any vertex is $2\pi$, we have $\alpha=\frac{2\pi}{3}$, $\delta<\alpha$, $\gamma>\alpha$ and $\beta<\alpha$, and we get $\beta<\gamma$ and $\delta<\epsilon$. If $\epsilon<\alpha$, then we get the similar but reversed inequalities. Since this always contradicts the geometrical constraint in \cite[Lemma 3]{ccy1}, we conclude that $\epsilon=\alpha$, and furthermore all angles are equal to $\alpha$. Therefore the tiling is the regular dodecahedron. 

The type 4 tiling is given on the right of Figure \ref{T14}. The normal edges have the same length, and the thick edges have the same length. We have $\beta+2\epsilon=\gamma+2\delta=2\pi$, which implies $\beta<\gamma$ if and only if $\epsilon>\delta$, and $\beta>\gamma$ if and only if $\epsilon<\delta$. Then \cite[Lemma 3]{ccy1} implies $\beta=\gamma$ and $\delta=\epsilon$. This further implies all angles are equal, so that the tiling is again the regular dodecahedron. 

\bigskip

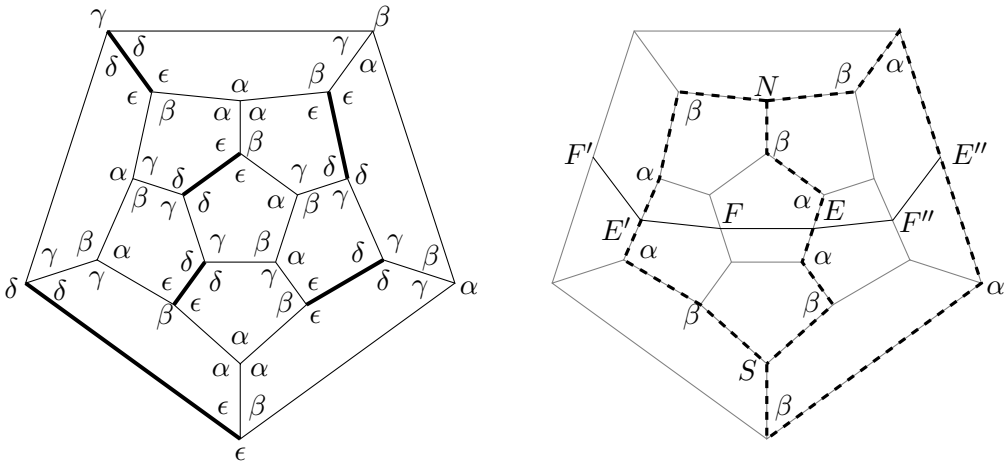
\begin{figure}[htp]
\centering
\begin{tikzpicture}[scale=1]

\foreach \a in {0,1}
\foreach \x in {1,...,5}
{
\begin{scope}[xshift=7*\a cm]
\coordinate (X\x A\a) at (18+72*\x:0.5);
\coordinate (Y\x A\a) at (6+72*\x:0.95);
\coordinate (Z\x A\a) at (-38+72*\x:0.95);
\coordinate (U\x A\a) at (8+72*\x:1.35);
\coordinate (V\x A\a) at (-44+72*\x:1.35);
\coordinate (W\x A\a) at (-18+72*\x:1.65);
\coordinate (P\x A\a) at (18+72*\x:1.7);
\coordinate (Q\x A\a) at (-11+72*\x:2.1);
\coordinate (R\x A\a) at (47+72*\x:2.1);
\coordinate (S\x A\a) at (-13+72*\x:2.6);
\coordinate (T\x A\a) at (49+72*\x:2.6);
\coordinate (O\x A\a) at (-18+72*\x:3.2);
\coordinate (P\a) at (-60:3.5);
\end{scope}
}

\foreach \x in {1,...,5}
\draw[rotate=72*\x]
	(90:0.8) -- (18:0.8) -- (18:1.5) -- (54:2) -- (90:1.5)
	(54:2) -- (54:3) -- (126:3);
	
\draw[line width=1.5]
	(90:0.8) -- (162:0.8)
	(18:1.5) -- (54:2)
	(234:0.8) -- (234:1.5)
	(-54:1.5) -- (-18:2)
	(126:2) -- (126:3)
	(198:3) -- (270:3);


\node at (X1A0) {\small $\epsilon$};
\node at (X2A0) {\small $\delta$};
\node at (X3A0) {\small $\gamma$};
\node at (X4A0) {\small $\beta$};
\node at (X5A0) {\small $\alpha$};

\node at (Y1A0) {\small $\beta$};
\node at (Z1A0) {\small $\gamma$};
\node at (U1A0) {\small $\alpha$};
\node at (V1A0) {\small $\delta$};
\node at (W1A0) {\small $\epsilon$};

\node at (Y2A0) {\small $\delta$};
\node at (Z2A0) {\small $\epsilon$};
\node at (U2A0) {\small $\gamma$};
\node at (V2A0) {\small $\alpha$};
\node at (W2A0) {\small $\beta$};

\node at (Y3A0) {\small $\delta$};
\node at (Z3A0) {\small $\gamma$};
\node at (U3A0) {\small $\epsilon$};
\node at (V3A0) {\small $\beta$};
\node at (W3A0) {\small $\alpha$};

\node at (Y4A0) {\small $\gamma$};
\node at (Z4A0) {\small $\delta$};
\node at (U4A0) {\small $\beta$};
\node at (V4A0) {\small $\epsilon$};
\node at (W4A0) {\small $\alpha$};

\node at (Y5A0) {\small $\beta$};
\node at (Z5A0) {\small $\alpha$};
\node at (U5A0) {\small $\gamma$};
\node at (V5A0) {\small $\epsilon$};
\node at (W5A0) {\small $\delta$};

\node at (P1A0) {\small $\alpha$};
\node at (Q1A0) {\small $\beta$};
\node at (R1A0) {\small $\epsilon$};
\node at (S1A0) {\small $\gamma$};
\node at (T1A0) {\small $\delta$};

\node at (P2A0) {\small $\alpha$};
\node at (Q2A0) {\small $\epsilon$};
\node at (R2A0) {\small $\beta$};
\node at (S2A0) {\small $\delta$};
\node at (T2A0) {\small $\gamma$};

\node at (P3A0) {\small $\beta$};
\node at (Q3A0) {\small $\gamma$};
\node at (R3A0) {\small $\alpha$};
\node at (S3A0) {\small $\delta$};
\node at (T3A0) {\small $\epsilon$};

\node at (P4A0) {\small $\epsilon$};
\node at (Q4A0) {\small $\alpha$};
\node at (R4A0) {\small $\delta$};
\node at (S4A0) {\small $\beta$};
\node at (T4A0) {\small $\gamma$};

\node at (P5A0) {\small $\delta$};
\node at (Q5A0) {\small $\gamma$};
\node at (R5A0) {\small $\epsilon$};
\node at (S5A0) {\small $\beta$};
\node at (T5A0) {\small $\alpha$};

\node at (O1A0) {\small $\beta$};
\node at (O2A0) {\small $\gamma$};
\node at (O3A0) {\small $\delta$};
\node at (O4A0) {\small $\epsilon$};
\node at (O5A0) {\small $\alpha$};


\begin{scope}[xshift=7cm]

\foreach \x in {1,...,5}
\draw[gray,rotate=72*\x]
	(90:0.8) -- (18:0.8) -- (18:1.5) -- (54:2) -- (90:1.5)
	(54:2) -- (54:3) -- (126:3);

\draw[dashed,very thick]
	(90:1.5) -- (90:0.8) -- (18:0.8) -- (-54:0.8) -- (-54:1.5) -- (-90:2)
	(90:1.5) -- (126:2) -- (162:1.5) -- (198:2) -- (234:1.5) -- (-90:2)
	(90:1.5) -- (54:2) -- (54:3) -- (-18:3) -- (-90:3) -- (-90:2);

\draw
	(18:2.43) -- (-3:1.68) -- (-18:0.65) -- (198:0.65) -- (183:1.68) -- (162:2.43);
	
\node at (0.9,0.05) {\small $E$};
\node at (-0.45,0) {\small $F$};
\node at (2,-0.1) {\small $F''$};
\node at (-2,-0.2) {\small $E'$};
\node at (-2.5,0.8) {\small $F'$};
\node at (2.7,0.8) {\small $E''$};

\end{scope}

\node at (P1A1) {\small $N$};
\node at (R3A1) {\small $S$};

\node at (Y1A1) {\small $\beta$};
\node at (X5A1) {\small $\alpha$};
\node at (Z5A1) {\small $\alpha$};
\node at (U4A1) {\small $\beta$};

\node at (W2A1) {\small $\beta$};
\node at (P2A1) {\small $\alpha$};
\node at (W3A1) {\small $\alpha$};
\node at (P3A1) {\small $\beta$};

\node at (Q1A1) {\small $\beta$};
\node at (T5A1) {\small $\alpha$};
\node at (O5A1) {\small $\alpha$};
\node at (S4A1) {\small $\beta$};

\end{tikzpicture}
\caption{Type 2 tiling.}
\label{T2}
\end{figure}

\noindent{\bf 2.} 
The simple argument for types 1 and 4 does not work for the types 2 and type 3. The left of Figure \ref{T2} is the type 2 tiling, again with normal edges having the same length, and the thick edges having the same length. We denote the two $\alpha^3$-vertices by $N$ and $S$ (north and south poles). On the right of Figure \ref{T2}, the three dashed paths connecting $N$ and $S$ are congruent, and can be moved to each other by rotating by the angle $\frac{2\pi}{3}$ or $\frac{4\pi}{3}$. Therefore $N$ and $S$ are antipodal points. Moreover, the dashed paths are also congruent with respect to the exchange of $N$ and $S$. This implies that the equator of the sphere (with respect to the two poles) cuts three indicated edges at the middle points $E,E',E''$, which are of distance $\frac{1}{3}2\pi=\frac{2\pi}{3}$ from each other. Moreover, the distance between $E$ and $N$ is $\frac{\pi}{2}$. Moreover, since all tiles are congruent, we also find that the equator cuts the other three indicated edges at the middle points $F,F',F''$, and $F$s are the middle points between the $E$s, so that the distance between $E$ and $F$ is $\frac{1}{2}\frac{2\pi}{3}=\frac{\pi}{3}$.

\begin{figure}[htp]
\centering
\begin{tikzpicture}


\draw
	(-1.2,1.2) node[above=-2] {\small $N$} 
	-- node[fill=white,inner sep=1] {\small $a$} 
	(0,0) -- node[fill=white,inner sep=1] {\small $a$}
	(1.2,1.2) -- node[fill=white,inner sep=1] {\small $\frac{a}{2}$}
	(2.2,0.6) node[right=-2] {\small $E$};
\draw[dashed]
	(-1.2,1.2) -- (2.2,0.6);
	
\node at (0,0.3) {\small $\beta$};
\node at (1.2,0.95) {\small $\alpha$};
\node[fill=white,inner sep=1] at (0,1) {\small $\frac{\pi}{2}$};


\begin{scope}[xshift=4cm]

\draw
	(-0.5,1.2) node[above=-2] {\small $F$} 
	-- node[fill=white,inner sep=1] {\small $\frac{a}{2}$} 
	(0,0) -- node[fill=white,inner sep=1] {\small $a$}
	(2,0) -- node[fill=white,inner sep=1] {\small $\frac{a}{2}$}
	(2.5,1.2) node[above=-2] {\small $E$};
\draw[dashed]
	(-0.5,1.2) -- node[fill=white,inner sep=1] {\small $\frac{\pi}{3}$} 
	(2.5,1.2);
	
\node at (0.1,0.2) {\small $\gamma$};
\node at (1.9,0.2) {\small $\beta$};

\end{scope}


\begin{scope}[xshift=8cm]

\draw
	(-0.5,1.5) -- node[fill=white,inner sep=1] {\small $a$}
	(0,0) -- node[fill=white,inner sep=1] {\small $b$}
	(1.5,0) -- node[fill=white,inner sep=1] {\small $c$}
	(2.5,1);
\draw[dashed]
	(2.5,1) -- node[fill=white,inner sep=1] {\small $x$} 
	(-0.5,1.5);
	
\node at (0.2,0.2) {\small $\phi$};
\node at (1.4,0.2) {\small $\psi$};

\end{scope}

\end{tikzpicture}
\caption{The length of the fourth edge in a quadrilateral.}
\label{quad}
\end{figure}
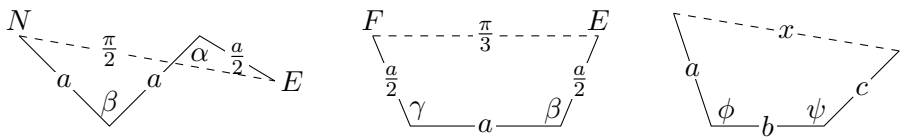

Let $a$ be the length of the normal edges. The points $N$ and $E$ as related by three edges on the left of Figure \ref{quad}, and the points $E$ and $F$ are related by another three edges in the middle of Figure \ref{quad}. The distances between the end points give two equalities relating $a$ and $\beta$. Our idea is to solve the two equalities to get the specific values of $a$ and $\beta$. 

On the right of Figure \ref{quad} is a quadrilateral with three consecutive edges $a,b,c$ and angles $\phi,\psi$ between these edges. The spherical trigonometry gives the fourth edge $x$ by
\begin{align}
\cos x
&=\cos a\cos b\cos c+\sin b(\sin a\cos c\cos\phi+\cos a\sin c\cos\psi) \nonumber \\
&\quad +\sin a\sin c(\sin\phi\sin\psi-\cos b\cos \phi\cos\psi). \label{quadlaw}
\end{align}
Then we may apply the formula to the left of Figure \ref{quad} with
\[
x=\frac{\pi}{2},\;
\phi=\beta,\;
\psi=2\pi-\alpha=\frac{4\pi}{3}.
\]
To make the equality more robust for numerical calculation, we may further divide $\cos\frac{a}{2}$ (which is nonzero because $0<a<\pi$ in the tiling) to get
\begin{align*}
0
&=\cos^2a+2\sin \frac{a}{2}\left(\sin a\cos \frac{a}{2}\cos\beta+\cos a\sin \frac{a}{2}\cos\frac{4\pi}{3}\right) \nonumber \\
&\quad +2\sin^2\frac{a}{2}\left(\sin\beta\sin\frac{4\pi}{3}-\cos a\cos \beta\cos\frac{4\pi}{3}\right).
\end{align*}
Similarly, we may apply \eqref{quadlaw} to the middle of Figure \ref{quad} with
\[
x=\frac{\pi}{3},\;
\phi=\gamma=\frac{4\pi}{3}-\beta,\;
\psi=\beta.
\]
The equality we get is already robust for numerical calculation.

\begin{figure}[htp]
\centering
\begin{tikzpicture}

\begin{scope}[scale=1]

\foreach \x in {1,...,5}
{
\draw[rotate=72*\x]
	(90:0.8) -- (18:0.8) -- (18:1.5) -- (54:2) -- (90:1.5)
	(54:2) -- (54:3) -- (126:3);

\coordinate (X\x A1) at (18+72*\x:0.5);
\coordinate (Y\x A1) at (6+72*\x:0.95);
\coordinate (Z\x A1) at (-38+72*\x:0.95);
\coordinate (U\x A1) at (8+72*\x:1.35);
\coordinate (V\x A1) at (-44+72*\x:1.35);
\coordinate (W\x A1) at (-18+72*\x:1.65);
\coordinate (P\x A1) at (18+72*\x:1.7);
\coordinate (Q\x A1) at (-11+72*\x:2.1);
\coordinate (R\x A1) at (47+72*\x:2.1);
\coordinate (S\x A1) at (-13+72*\x:2.6);
\coordinate (T\x A1) at (49+72*\x:2.6);
\coordinate (O\x A1) at (-18+72*\x:3.2);
\coordinate (P1) at (-60:3.5);
}

\draw[line width=1.5]
	(-54:0.8) -- (234:0.8)
	(18:0.8) -- (18:1.5)
	(162:0.8) -- (162:1.5)
	(-18:2) -- (-18:3)
	(198:2) -- (198:3)
	(54:3) -- (126:3);


\node at (X1A1) {\small $\alpha$};
\node at (X2A1) {\small $\beta$};
\node at (X3A1) {\small $\delta$};
\node at (X4A1) {\small $\epsilon$};
\node at (X5A1) {\small $\gamma$};

\node at (Y1A1) {\small $\gamma$};
\node at (Z1A1) {\small $\epsilon$};
\node at (U1A1) {\small $\alpha$};
\node at (V1A1) {\small $\delta$};
\node at (W1A1) {\small $\beta$};

\node at (Y2A1) {\small $\delta$};
\node at (Z2A1) {\small $\beta$};
\node at (U2A1) {\small $\epsilon$};
\node at (V2A1) {\small $\alpha$};
\node at (W2A1) {\small $\gamma$};

\node at (Y3A1) {\small $\beta$};
\node at (Z3A1) {\small $\delta$};
\node at (U3A1) {\small $\alpha$};
\node at (V3A1) {\small $\epsilon$};
\node at (W3A1) {\small $\gamma$};

\node at (Y4A1) {\small $\epsilon$};
\node at (Z4A1) {\small $\delta$};
\node at (U4A1) {\small $\gamma$};
\node at (V4A1) {\small $\beta$};
\node at (W4A1) {\small $\alpha$};

\node at (Y5A1) {\small $\epsilon$};
\node at (Z5A1) {\small $\gamma$};
\node at (U5A1) {\small $\delta$};
\node at (V5A1) {\small $\alpha$};
\node at (W5A1) {\small $\beta$};

\node at (P1A1) {\small $\alpha$};
\node at (Q1A1) {\small $\gamma$};
\node at (R1A1) {\small $\beta$};
\node at (S1A1) {\small $\epsilon$};
\node at (T1A1) {\small $\delta$};

\node at (P2A1) {\small $\gamma$};
\node at (Q2A1) {\small $\alpha$};
\node at (R2A1) {\small $\epsilon$};
\node at (S2A1) {\small $\beta$};
\node at (T2A1) {\small $\delta$};

\node at (P3A1) {\small $\gamma$};
\node at (Q3A1) {\small $\epsilon$};
\node at (R3A1) {\small $\alpha$};
\node at (S3A1) {\small $\delta$};
\node at (T3A1) {\small $\beta$};

\node at (P4A1) {\small $\beta$};
\node at (Q4A1) {\small $\alpha$};
\node at (R4A1) {\small $\delta$};
\node at (S4A1) {\small $\gamma$};
\node at (T4A1) {\small $\epsilon$};

\node at (P5A1) {\small $\beta$};
\node at (Q5A1) {\small $\delta$};
\node at (R5A1) {\small $\alpha$};
\node at (S5A1) {\small $\epsilon$};
\node at (T5A1) {\small $\gamma$};

\node at (O1A1) {\small $\epsilon$};
\node at (O2A1) {\small $\delta$};
\node at (O3A1) {\small $\beta$};
\node at (O4A1) {\small $\alpha$};
\node at (O5A1) {\small $\gamma$};

\end{scope}


\begin{scope}[shift={(6cm,0.5cm)}]

\draw
	(-1.2,1.2) node[above=-2] {\small $N$} 
	-- node[fill=white,inner sep=1] {\small $a$} 
	(0,0) -- node[fill=white,inner sep=1] {\small $a$}
	(1.2,1.2) -- node[fill=white,inner sep=1] {\small $\frac{a}{2}$}
	(2.2,0.6) node[right=-2] {\small $E$};
\draw[dashed]
	(-1.2,1.2) -- (2.2,0.6);
	
\node at (0,0.3) {\small $\gamma$};
\node at (1.2,0.95) {\small $\gamma$};
\node[fill=white,inner sep=1] at (0,1) {\small $\frac{\pi}{2}$};

\end{scope}


\begin{scope}[shift={(6cm,-2cm)}]

\draw
	(-1.2,1.2) node[above=-2] {\small $N$} 
	-- node[fill=white,inner sep=1] {\small $a$} 
	(0,0) -- node[fill=white,inner sep=1] {\small $a$}
	(1.2,1.2) -- node[fill=white,inner sep=1] {\small $\frac{a}{2}$}
	(2.2,0.6) node[right=-2] {\small $E$};
\draw[dashed]
	(-1.2,1.2) -- (2.2,0.6);
	
\node at (0,0.35) {\small $\beta$};
\node at (1.2,0.95) {\small $\beta$};
\node[fill=white,inner sep=1] at (0,1) {\small $\frac{\pi}{2}$};
\end{scope}

\end{tikzpicture}
\caption{Type 3 tiling.}
\label{T3}
\end{figure}
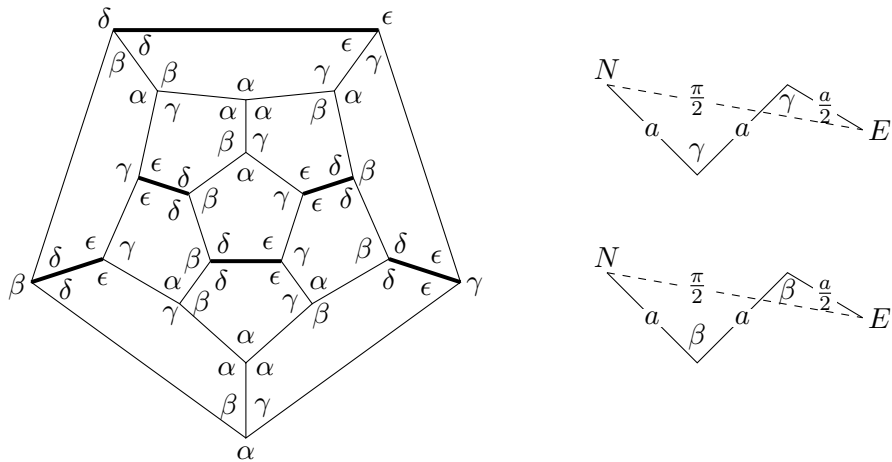

Now we turn to the type 3 tiling, given on the left of Figure \ref{T3}. Similar to the type 2 tiling, we find three congruent paths connecting the two poles, with turning angles $\gamma,\gamma,\gamma,\gamma$ instead of $\beta,\alpha,\alpha,\beta$. By the symmetry of exchanging $(\beta,\delta)$ with $(\gamma,\epsilon)$, we find another three congruent paths connecting the two poles, with turning angles $\beta,\beta,\beta,\beta$. We also have the symmetry of exchanging the two poles. Therefore we conclude that the equator of the sphere again cuts six edges at the middle points $E,E',E'',F,F',F''$. This gives two pictures on the right of Figure \ref{T3}. So we may substitute the following two sets of data into \eqref{quadlaw}
\begin{align*}
x&=\frac{\pi}{2},\;
\phi=\gamma=\frac{4\pi}{3}-\beta,\;
\psi=2\pi-\gamma=\frac{2\pi}{3}+\beta;\\
x&=\frac{\pi}{2},\;
\phi=\beta,\;
\psi=2\pi-\beta.
\end{align*}
We may further divide $\cos\frac{a}{2}$ to get more robust equalities for the numerical calculation.

\bigskip

\begin{figure}[htp]
\centering
\begin{tikzpicture}[scale=0.4]

\pgftext{\includegraphics{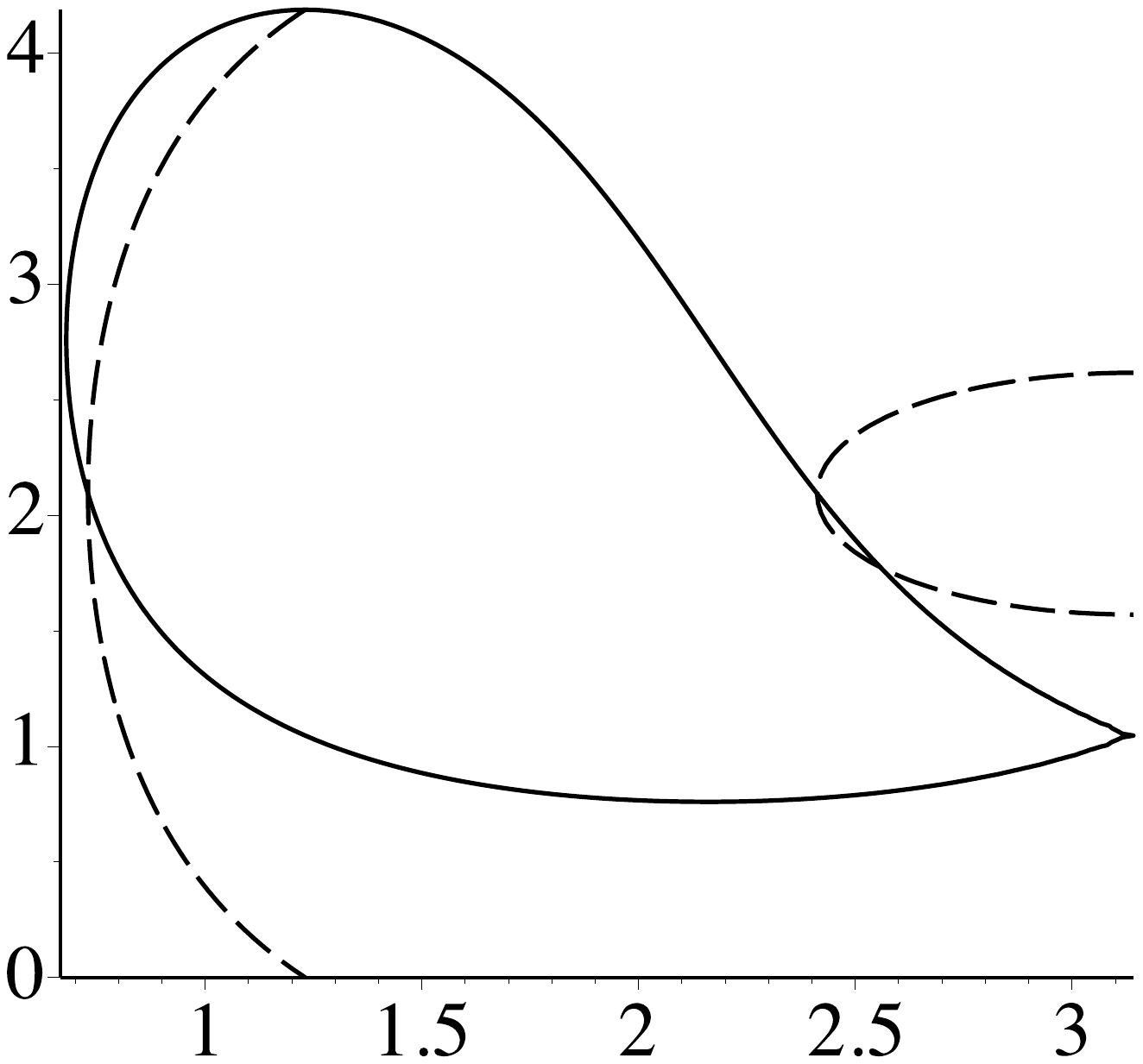}}

\node at (-6.1,6.7) {\small $\beta$};
\node at (7,-5.2) {\small $a$};

\node at (-5.2,0.6) {\small $P$};
\node at (-3.3,6.8) {\small $Q$};
\node at (2.3,0.3) {\small $R$};
\node at (3.3,-0.9) {\small $S$};

\begin{scope}[xshift=18cm]

\pgftext{\includegraphics{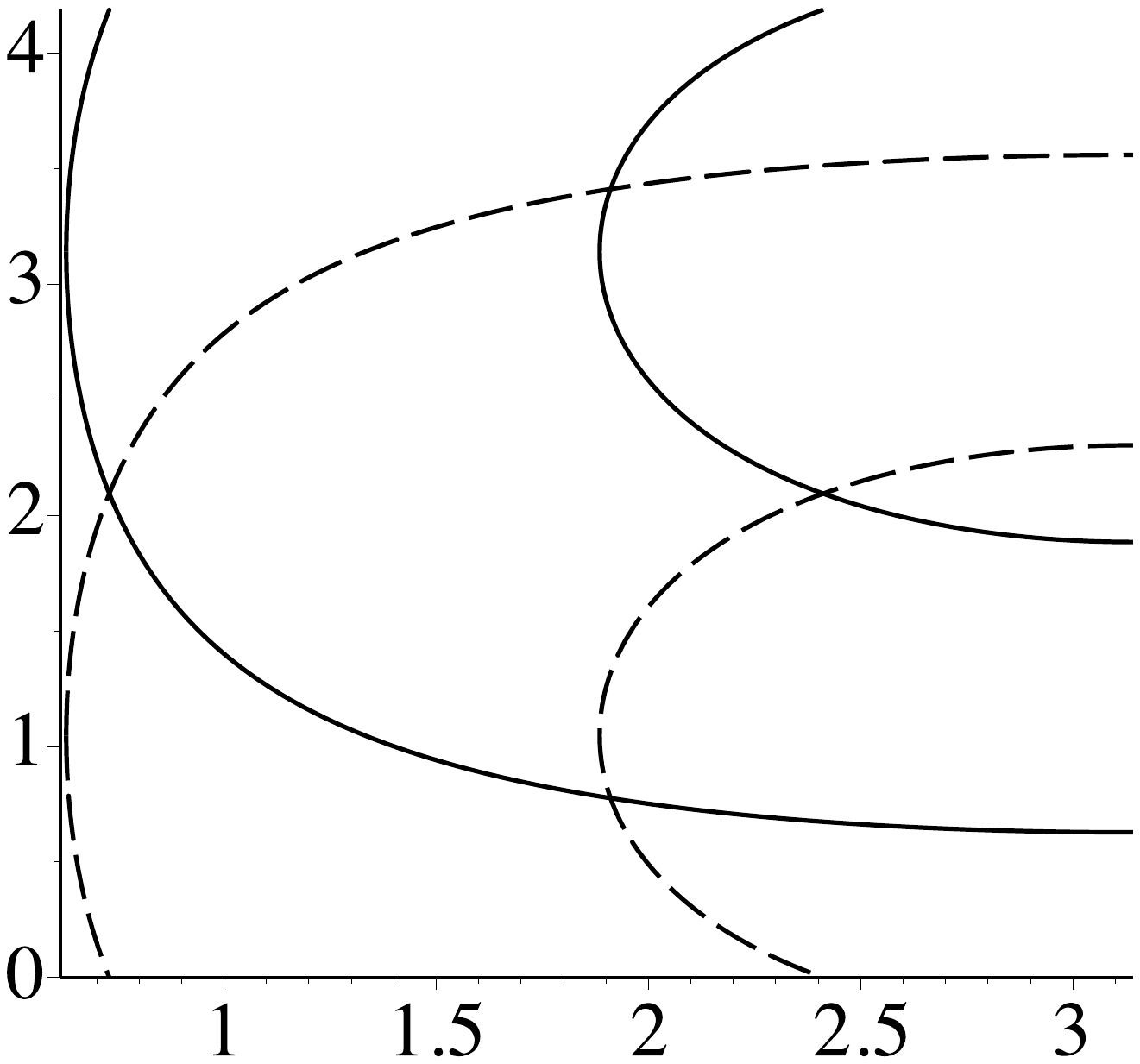}}

\node at (-6,6.7) {\small $\beta$};
\node at (7.1,-5.2) {\small $a$};

\node at (-4.9,0.6) {\small $P$};
\node at (0.1,4.6) {\small $Q$};
\node at (2.9,1.1) {\small $R$};
\node at (0.1,-3.5) {\small $S$};

\end{scope}

\end{tikzpicture}
\caption{$a$ and $\beta$ for type 2 and type 3 tilings.}
\label{numerical}
\end{figure}

\noindent{\bf 3.} 
Figure \ref{numerical} gives the results of the numerical calculations for the type 2 and 3 tilings. For each type, the two equations for $(a,\beta)\in (0,\pi)\times(0,\frac{4\pi}{3})$ give two curves, and their intersections are the solutions. On the left, we get four solutions $P,Q,R,S$ for the type 2 tiling. The solution $P$ corresponds to the regular dodecahedron. The solution $Q$ has $\beta=\frac{4\pi}{3}$, which is not allowed. Moreover, we note that if $\frac{\pi}{3}<\beta<\frac{4\pi}{3}$, then $\gamma=\frac{4\pi}{3}-\beta<\pi$ and $\delta=\frac{\pi}{3}+\frac{\beta}{2}<\pi$. This implies that the isosceles triangle with two sides $a$ and angle $\alpha=\frac{2\pi}{3}$ between the two sides is contained in the pentagonal tile of area $\frac{4\pi}{12}=\frac{\pi}{3}$. The geometrical fact means $\cos a>\frac{1}{3}$, or $a<1.231$. Since the solutions $R$ and $S$ satisfy $\frac{\pi}{3}<\beta<\frac{4\pi}{3}$ and $a\ge 1.231$, they are not geometrically realistic (these are pentagons with boundary edges crossing each other). We conclude that the regular dodecahedron is the only geometrically realistic solution.

On the right of Figure \ref{numerical}, we get four solutions $P,Q,R,S$ for the type 3 tiling. The solution $P$ corresponds to the regular dodecahedron. Moreover, it is easy to see that $\delta<\pi$ and $\epsilon< \pi$. As we argued for the type 2 tiling, this implies $\cos a>\frac{1}{3}$, or $a<1.231$. Therefore $Q,R,S$ are not geometrically realistic solutions. We conclude again that the regular dodecahedron is the only geometrically realistic solution.

\bigskip

\begin{figure}[htp]
\centering
\begin{tikzpicture}[scale=1]

\foreach \x in {1,...,5}
{
\coordinate (A\x) at (18+72*\x:0.8);
\coordinate (B\x) at (18+72*\x:1.5);
\coordinate (C\x) at (-18+72*\x:2);
\coordinate (D\x) at (-18+72*\x:3);

\coordinate (X\x) at (18+72*\x:0.5);
\coordinate (Y\x) at (6+72*\x:0.95);
\coordinate (Z\x) at (-38+72*\x:0.95);
\coordinate (U\x) at (8+72*\x:1.35);
\coordinate (V\x) at (-44+72*\x:1.35);
\coordinate (W\x) at (-18+72*\x:1.7);
\coordinate (P\x) at (18+72*\x:1.75);
\coordinate (Q\x) at (-11+72*\x:2.1);
\coordinate (R\x) at (47+72*\x:2.1);
\coordinate (S\x) at (-13+72*\x:2.6);
\coordinate (T\x) at (49+72*\x:2.6);
\coordinate (O\x) at (-18+72*\x:3.3);
}

\draw
	(A1) -- (A2) -- (A3)
	(A2) -- (B2)
	(B1) -- (C1) -- (B5)
	(C1) -- (D1)
	(C4) -- (B4) -- (C5)
	(A4) -- (B4)
	(D2) -- (D3) -- (D4)
	(C3) -- (D3);

\draw[line width=1.5]
	(A1) -- (A5) -- (A4)
	(A5) -- (B5)
	(B1) -- (C2) -- (B2)
	(C2) -- (D2)
	(C3) -- (B3) -- (C4)
	(A3) -- (B3)
	(D1) -- (D5) -- (D4)
	(C5) -- (D5);

\draw[dashed]
	(A3) -- (A4)
	(A1) -- (B1)
	(B2) -- (C3)
	(B5) -- (C5)
	(C4) -- (D4)
	(D1) -- (D2);

\fill
	(A2) circle (0.1)
	(A5) circle (0.1)
	(B3) circle (0.1)
	(B4) circle (0.1)
	(C1) circle (0.1)
	(C2) circle (0.1)
	(D3) circle (0.1)
	(D5) circle (0.1);

\node at (X1) {\small $\delta$};
\node at (X2) {\small $\alpha$};
\node at (X3) {\small $\beta$};
\node at (X4) {\small $\gamma$};
\node at (X5) {\small $\alpha$};

\node at (Y1) {\small $\gamma$};
\node at (Z1) {\small $\alpha$};
\node at (U1) {\small $\beta$};
\node at (V1) {\small $\delta$};
\node at (W1) {\small $\alpha$};

\node at (Y2) {\small $\alpha$};
\node at (Z2) {\small $\beta$};
\node at (U2) {\small $\delta$};
\node at (V2) {\small $\gamma$};
\node at (W2) {\small $\alpha$};

\node at (Y3) {\small $\delta$};
\node at (Z3) {\small $\alpha$};
\node at (U3) {\small $\alpha$};
\node at (V3) {\small $\beta$};
\node at (W3) {\small $\gamma$};

\node at (Y4) {\small $\beta$};
\node at (Z4) {\small $\gamma$};
\node at (U4) {\small $\alpha$};
\node at (V4) {\small $\alpha$};
\node at (W4) {\small $\delta$};

\node at (Y5) {\small $\alpha$};
\node at (Z5) {\small $\delta$};
\node at (U5) {\small $\gamma$};
\node at (V5) {\small $\alpha$};
\node at (W5) {\small $\beta$};

\node at (P1) {\small $\delta$};
\node at (Q1) {\small $\alpha$};
\node at (R1) {\small $\alpha$};
\node at (S1) {\small $\beta$};
\node at (T1) {\small $\gamma$};

\node at (P2) {\small $\gamma$};
\node at (Q2) {\small $\alpha$};
\node at (R2) {\small $\beta$};
\node at (S2) {\small $\delta$};
\node at (T2) {\small $\alpha$};

\node at (P3) {\small $\alpha$};
\node at (Q3) {\small $\delta$};
\node at (R3) {\small $\gamma$};
\node at (S3) {\small $\alpha$};
\node at (T3) {\small $\beta$};

\node at (P4) {\small $\alpha$};
\node at (Q4) {\small $\beta$};
\node at (R4) {\small $\delta$};
\node at (S4) {\small $\gamma$};
\node at (T4) {\small $\alpha$};

\node at (P5) {\small $\beta$};
\node at (Q5) {\small $\gamma$};
\node at (R5) {\small $\alpha$};
\node at (S5) {\small $\alpha$};
\node at (T5) {\small $\delta$};

\node at (O1) {\small $\gamma$};
\node at (O2) {\small $\beta$};
\node at (O3) {\small $\alpha$};
\node at (O4) {\small $\delta$};
\node at (O5) {\small $\alpha$};

\end{tikzpicture}
\caption{Type 5 tiling.}
\label{T5}
\end{figure}
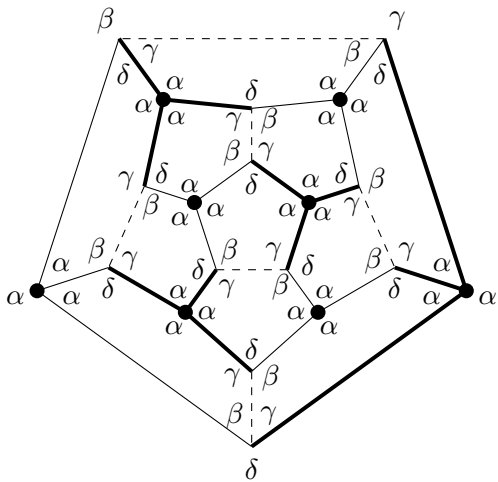

\noindent{\bf 4.} 
By \cite{gsy} and the argument above, the only deformed dodecahedron tiling is of type 5, given in Figure \ref{T5}. We observe that the tiling is symmetric with respect to the three fold rotations at the eight dotted vertices. By combining these rotations, we can move any tile to any other tile. Therefore the tiling is {\em isohedral}, in the sense that the symmetry group acts transitively on the tiles.

The full symmetry of the deformed dodecahedron tiling depends on the lengths $a,b,c$ of the normal edges, thick edges, and dashed edges. The symmetry group $G$ may be obtained by observing that the transitivity of the action implies $G/H=\{12\text{ tiles}\}$, where $H$ is the subgroup of symmetries that preserve one tile (say the center tile).

For the case that $a,b,c$ are distinct, the only symmetry of the center tile is the identity, and fixing the center tile implies fixing all the tiles. Therefore $H$ is the trivial group, and the order of $G$ is $12$. It turns out that $G$ is the chiral tetrahedral group $T$. 

For the case $a=c\ne b$, we know $\alpha\ne\beta$ by \cite[Lemma 3]{ccy1}. This implies that $H$ is still the trivial group, so that $G=T$. The same happens for the case $b=c\ne a$.

The case $a=b\ne c$ is the pyritohedron. We have $\beta=\gamma$ by \cite[Lemma 3]{ccy1}, so that the tile is symmetric with respect to the flipping that preserves the $\delta$ angle. Moreover, the flipping of the center tile determines the action on all the other tiles. Therefore $H$ has order $2$, and $G$ has order $24$. In fact, $G$ is the pyritohedral group $T_h$. 

The case $a=b=c$ is the regular dodecahedron. The subgroup $H$ is the symmetry group of the regular pentagon, which is the dihedral group of order $10$. Moreover, $G$ is the icosahedral group $I_h$.

The symmetry of deformed dodecahedron tiling is summarized in Table~\ref{symmetry}.

\begin{table}[htp]
\centering
\begin{tabular}{|c|c|c|}
\hline 
edges & symmetry & order \\
\hline   
$a\ne b$ 
& $T$ 
& $12$ \\
\hline 
$a=b\ne c$
& $T_h$ 
& $24$ \\
\hline 
$a=b=c$
& $I_h$ 
& $120$ \\
\hline 
\end{tabular}
\caption{Symmetry of deformed dodecahedron tiling.}
\label{symmetry}
\end{table}

\bigskip

\noindent{\bf 5.} 
The deformed dodecahedron tiling may be compared with the other deformed platonic solids. By \cite[Theorem 1]{ua}, the deformed tetrahedron $\bullet F_4$ allows $2$ free variables and has the dihedral group $D_2$ as the symmetry group. (The symmetry may be bigger if some edges become equal.) The deformed octahedron $\bullet G_8$ allows $1$ free variable and has the dihedral group $D_{2d}$ as the symmetry group. The icosahedron $H_{20}$ is not deformable (i.e., rigid) and has the icosahedral group $I_h$ as the symmetry group. By \cite[Theorem 2]{Akama3} and \cite[Theorem 1]{Akama2}, the deformed cube tiling allows $2$ free variables and has the dihedral group $D_3$ as the symmetry group.

\begin{figure}[htp]
\centering
\begin{tikzpicture}[scale=1]


\coordinate (A) at (0,0);
\coordinate (B) at (1,1.6);
\coordinate (C) at (1.8,-0.8);
\coordinate (D) at (2.4,0.4);

\draw
	(A) -- (D);

\fill[white]
	(1.4,0.2) circle (0.1);

\draw
	(B) -- (C);
		
\draw[dashed]
	(A) -- (B)
	(C) -- (D);
	
\draw[line width=1.5]
	(A) -- (C)
	(B) -- (D);

\draw[densely dotted]
	(1.2,0.2) -- (1.4,0.4) -- ++(0.7,0.7) node[above right=-3] {\small $2$}
	(1.7,1) -- (0.9,-0.4) -- ++(-0.2,-0.35) node[below left=-3] {\small $2$}
	(2.1,-0.2) -- (0.5,0.8) -- ++ (-0.4,0.25) node[above left=-3] {\small $2$}
	;
	
\node at (1.2,-1.5) {$\bullet F_4$};


\begin{scope}[shift={(5cm,0.2cm)}]

\coordinate (A1) at (0,1.2);
\coordinate (A2) at (0,-1.2);
\coordinate (B1) at (-1,-0.6);
\coordinate (D1) at (1,-0.6);
\coordinate (B2) at (1,0.6);
\coordinate (D2) at (-1,0.6);
\coordinate (C1) at (0.1,0.3);
\coordinate (C2) at (-0.1,-0.3);

\draw
	(A1) -- (B2)
	(A1) -- (C1)
	(A1) -- (D2)
	(A2) -- (B1)
	(A2) -- (C2)
	(A2) -- (D1);

\draw[dashed]
	(C2) -- (D2);

\draw[line width=1.5]
	(B2) -- (C2);

\fill[white]
	(0.33,0.07) circle (0.1)
	(-0.33,0) circle (0.15);
	
\draw[dashed]
	(B1) -- (C1) 
	(D1) -- (B2);

\draw[line width=1.5]
	(C1) -- (D1)
	(D2) -- (B1);

\draw[densely dotted]
	(A2) -- (A1) -- ++(0,0.5) node[above=-3] {\small $3$}
	(-1,0) -- (1.3,0) node[right=-1] {\small $2$}
	(0.45,0.15) -- (-0.45,-0.15) -- ++(-0.9,-0.3) node[left=-2] {\small $2$}
	
	(0.55,-0.15) -- (-0.55,0.15) -- ++(-0.77,0.21) node[left=-2] {\small $2$};	

\node at (0,-1.7) {$\bullet$Cube};

\end{scope}


\begin{scope}[shift={(9cm,0.2cm)}]

\coordinate (A1) at (0,1.2);
\coordinate (A2) at (0,-1.2);
\coordinate (B1) at (-1.2,0.2);
\coordinate (B2) at (1.2,0.2);
\coordinate (C1) at (0.3,-0.5);
\coordinate (C2) at (-0.3,0.1);

\draw
	(B1) -- (A1) -- (B2)
	(C1) -- (A2) -- (C2);

\draw[dashed]
	(B2) -- (C2) -- (B1);
	
\fill[white]
	(0.2,0.1) circle (0.1)
	(-0.15,-0.25) circle (0.15);	

\draw[dashed]
	(B1) -- (C1) -- (B2);
	
\draw[line width=1.5]
	(A2) -- (B1) 
	(A2) -- (B2)
	(A1) -- (C1)
	(A1) -- (C2);
	
\draw[densely dotted]
	(-1.2,-1.2) rectangle (1.2,1.2)
	(-0.3,1.5) -- (0.3,0.9) -- (0.3,-1.5) -- (-0.3,-0.9) -- cycle
	(A2) -- (A1) -- ++(0,0.5) node[above=-3] {\small $2$}
	(0.45,0.15) -- (-0.45,-0.15) -- ++(-0.9,-0.3) node[left=-2] {\small $2$}
	(-0.75,0.15) -- (0.75,-0.15)-- ++(1,-0.2) node[right=-2] {\small $2$}
	;
	
\node at (-0.4,1.6) {\small m};
\node at (-1.3,1.3) {\small m};

\node at (0,-1.7) {$\bullet G_8$};

\end{scope}

\end{tikzpicture}
\caption{Deformed tetrahedron, cube, octahedron and their symmetries. The number $n$ attached to an axis means $n$-fold rotation. A plane with ``m'' means mirror reflection.}
\label{platonic}
\end{figure}
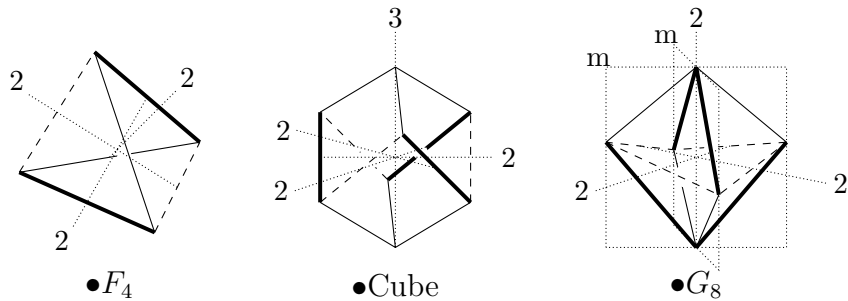

All the deformed platonic solids are isohedral. On the other hand, the trapezohedron of 12 faces, which is dual to an antiprism of 12 vertices, is isohedral when the tiles are convex, and may not be isohedral with concave tiles \cite[Theorem 5]{Akama3}.

Finally, we note that the deformed dodecahedron tiling is closely related to the dodecahedral shell in a so called Tsai-cluster~(or RTH-cluster) of quasicrystalline Cd-Yb alloy \cite{tgybt} and the dodecahedral cluster of water molecules deform~\cite{sx}. The isohedral property of the deformed dodecahedron tiling established in this note should shed some light on these chemical structures.

\end{document}